%Paper: alg-geom/9307003
%From: rf@shire.math.columbia.edu (Robert Friedman)
%Date: Wed, 14 Jul 93 12:08:53 -0400

%%Paper written in AMS-TeX
%%using the amsppt style

\input amstex

\define\scrO{\Cal O}
\define\Pee{{\Bbb P}}
\define\Zee{{\Bbb Z}}

\define\Id{\operatorname{Id}}
\define\Pic{\operatorname{Pic}}

\define\ad{\operatorname{ad}}

\define\proof{\demo{Proof}}
\define\endproof{\qed\enddemo}
\define\endstatement{\endproclaim}
\define\theorem#1{\proclaim{Theorem #1}}
\define\lemma#1{\proclaim{Lemma #1}}
\define\proposition#1{\proclaim{Proposition #1}}
\define\corollary#1{\proclaim{Corollary #1}}
\define\claim#1{\proclaim{Claim #1}}

\define\section#1{\specialhead #1 \endspecialhead}
\define\ssection#1{\medskip\noindent{\bf #1}}

\documentstyle{amsppt}
\NoBlackBoxes
\leftheadtext{Vector bundles  and $SO(3)$-invariants for elliptic surfaces II}
\rightheadtext{The case of even fiber degree}
\pageno=1
\topmatter
\title Vector bundles  and $SO(3)$-invariants for elliptic surfaces II:\\
The case of even fiber degree
\endtitle
\author Robert Friedman
\endauthor
\address Columbia University, New York, NY 10027
\endaddress
\email rf\@math.columbia.edu \endemail
\thanks Research partially supported by NSF grant  DMS-9203940
\endthanks
\endtopmatter

\document

\section{Introduction.}

Let $S$ be a simply connected elliptic surface with at most two multiple
fibers.
In this paper, the second in a series of three, we are concerned with
describing
moduli spaces of stable vector bundles $V$ over $S$ such that the restriction
of $c_1(V)$ to a general fiber has the smallest possible nonzero degree, namely
the product of the multiplicities, in the case where this product is even. We
then apply this study toward a partial calculation of the corresponding
Donaldson
polynomial invariants of $S$. Our goal is the completion of the $C^\infty$
classification of such surfaces, and the general outline of this classification
has been described in the introduction to Part I. Aside from quoting
a few results from Part I, this paper can however be read independently. On the
other hand, the methods of this paper draw heavily on the book [4], and many
arguments which are very similar to arguments in [4] are sketched or simply
omitted. Roughly speaking, the new ingredients in the proof consist of the
algebraic geometry of certain elliptic surfaces associated to $S$, which have a
single multiple fiber of multiplicity two and are birational to double covers
of rational ruled surfaces. The vector bundle parts of the argument run more or
less parallel to the arguments in [4], with a few new cases to analyze.

The outline of this paper is as follows.  In this
paper, we shall only be concerned with elliptic surfaces $S$ over $\Pee ^1$
with
multiple fibers of multiplicities $2m_1$ and $m_2$, where $m_2$ is odd, and
such
that there exists a divisor $\Delta$ on $S$ with $\Delta \cdot f = 2m_1m_2$,
the
minimum possible value, for a smooth fiber $f$. In this case, there is an
associated  surface $J^{m_1m_2}(S)$ defined in [3]. The surface $J^{m_1m_2}(S)$
fibers over $\Pee ^1$ and the fiber over a point $t$ lying under a smooth
fiber $f$ of $S$ is $J^{m_1m_2}(f)$, the set of line bundles  of degree
$m_1m_2$
on the fiber $f$ of $S$. The surface $J^{m_1m_2}(S)$ has an
involution defined by $\lambda \in J^{m_1m_2}(f) \mapsto \scrO _f(\Delta |f)
\otimes \lambda ^{-1}$. The quotient of  $J^{m_1m_2}(S)$ by this involution is
birational to a rational ruled surface $\Bbb F_N$, and we describe the geometry
of the double cover in detail. In Section 2, we describe the rough
classification of stable bundles $V$ on $S$ with $c_1(V)=\Delta$. To each such
bundle there is an associated bisection $C$ of $J^{m_1m_2}(S)$ which is
invariant under the involution, and so defines a section of the quotient ruled
surface. In Section 3, we show that for general bundles $V$, $V$ is determined
up to finite ambiguity by the section of the ruled surface and the choice of a
certain line bundle on the associated bisection $C$ of $J^{m_1m_2}(S)$. Thus, a
Zariski open subset  of the moduli  space fibers over an open subset of the
linear system of all sections on a certain rational ruled surface, and the
fibers are a number of copies of the Jacobian $J(C)$ of the bisection $C$ of
$J^{m_1m_2}(S)$. It is here that the asymmetry between $2m_1$ and $m_2$ becomes
apparent: the number of connected components of the fiber is just $m_2$. The
reasons for this are explained following Lemma 2.4. Finally, in Section 4 we
calculate the leading coefficient of a Donaldson polynomial invariant and show
that it contains an ``extra" factor of $m_2$.

\section{1. Geometry of elliptic surfaces.}

We fix notation for this paper. Let $\pi \: S\to \Pee ^1$ be an elliptic
surface
with two multiple fibers $F_{2m_1}$ and $F_{m_2}$ of multiplicities $2m_1$ and
$m_2$, where $\gcd (2m_1, m_2)=1$. We shall further assume that the reductions
of
the multiple fibers are smooth and that all other singular fibers are reduced
and
irreducible with just one ordinary double point. In other words, $S$ is nodal
in the terminology of [4]. Let $\kappa _S=\kappa \in H^2(S;\Zee)$ be the unique
class such that $2m_1m_2\kappa = [f]$, where $f$ is a general fiber of $\pi$.
Finally we shall assume that there is a $2m_1m_2$-section $\Delta$,
i.e\. a divisor $\Delta$, not necessarily effective, with $\Delta \cdot f =
2m_1m_2$,
or equivalently $\Delta \cdot \kappa = 1$. By [4], there always exist
such nodal elliptic surfaces. In particular $S$ is algebraic.

We shall be concerned with the associated elliptic surface $J^{m_1m_2}(S)$
defined in Section 2 of [3]. By  the discussion in Section 2 of [3],
$J^{m_1m_2}(S)$ has exactly one multiple fiber of multiplicity two, above the
point of $\Pee ^1$ corresponding to $F_{2m_1}$. We shall denote this fiber by
$F_2$. Moreover, given the $2m_1m_2$-section $\Delta$, there is an involution
of
$J^{m_1m_2}(S)$ defined on the generic fiber $S_\eta$ by $\lambda \mapsto
\scrO_{S_{\eta}}(\Delta)\otimes \lambda ^{-1}$. Since $J^{m_1m_2}(S)$ is
relatively minimal, this involution on the generic fiber extends to an
involution on $J^{m_1m_2}(S)$, which we shall denote by $\iota$.  The data of
$J^{m_1m_2}(S)$ and the involution $\iota$ do not depend on $\Delta$, but only
on the restriction of $\Delta$ to the generic fiber. A line bundle which has
the
same restriction as $\Delta$ to the generic fiber  differs from $\Delta$  by a
line bundle which is trivial on the generic fiber and thus is a multiple of
$\kappa$. Up to twisting by a line bundle which is divisible by 2, the only
possibilities are then $\Delta$ and $\Delta - \kappa$. Replacing $\Delta$ by
$\Delta - \kappa$ replaces $\Delta ^2$ by $\Delta ^2 -2$ and thus changes
$\Delta ^2 \mod 4$. Of course $\Delta ^2\equiv 1 \mod 2$ since
$$\Delta \cdot K_S = 2m_1m_2(p_g+1) - 2m_1-m_2\equiv 1 \mod 2.$$

Let us determine the fixed point set of $\iota$.

\lemma{1.1} The fixed point set of $\iota$ consists of a smooth $4$-section
together with two isolated fixed points on $F_2$.
\endstatement
\proof The fixed point set of an involution on the smooth surface
$J^{m_1m_2}(S)$
must consist of a smooth curve together with some isolated fixed points. For a
smooth nonmultiple fiber $f$, there are four divisors $\lambda$ such that
$2\lambda = f$.  Thus there is a component of the fixed point set which is a
4-section of $J^{m_1m_2}(S)$. Every other curve component of the fixed point
set has trivial restriction to the generic fiber and (since the curve component
is smooth) cannot meet the 4-section. Since all fibers are irreducible, there
can be no other component, and the remaining fixed points are isolated.

Let us consider the possibilities for isolated fixed points away from $F_2$.
In an analytic neighborhood of a nonmultiple fiber, there is a section of the
elliptic fibration. Using this section to make the local identification of
$J^{m_1m_2}(S)$ with $J^0(S)$, and since the group of local sections is
divisible, it is easy to see that, after a translation, we can assume that
$\iota$ corresponds to the involution $x\mapsto -x$. Direct inspection shows
that this involution has no isolated fixed points, even at the nodal fibers. To
handle the multiple fiber, the
explicit description of a neighborhood $X$ of a multiple fiber shows that after
a
base change of order two, say $\tilde X \to X$, we can assume that there is a
local section of $\tilde X$. The induced map from the central fiber of
$\tilde X$ to the central fiber of $X$ corresponds to taking the quotient by a
subgroup of order two. After a translation we can further assume that
the pulled back involution on $\tilde X$ is given by $x\mapsto -x$. Since
inverses commute with translations by a point of order two, the restriction of
$\iota$ to $F_2$ again has four fixed points. Two of these lie on the 4-section
(recall that a 4-section can  meet $F_2$ in at most
two distinct points) and the remaining two are isolated.
\endproof

Thus
there are two isolated fixed points of $\iota$ on $F_2$. If we blow these up
and
then take the quotient, the result will be a smooth surface $\Bbb F$ mapping to
$\Pee ^1$ whose fibers are smooth rational curves except over the point
corresponding to $F_2$ ({\smc Fig. 1}). Over this point, the fiber is a  curve
$\frak d_1+2\frak e+\frak d_2$, where $\frak d_1$ and $\frak d_2$ are the
images
of the exceptional curves, $\frak e$ is the image of $F_2$, and we have $(\frak
d_i)^2 = -2$ and $\frak e^2 = -1$, $\frak d_1\cdot \frak e =\frak  d_2\cdot
\frak e = 1$ and $\frak d_1\cdot  \frak d_2 = 0$. In particular we may
contract $\frak e$ and then either $\frak d_1$ or $\frak d_2$ to obtain a
rational ruled surface $\Bbb F_N$. We shall fix notation so that $\frak d_2$ is
the curve we contract and the resulting surface is $\Bbb F_N$. However, as we
shall see, it is important to keep in mind the symmetry between $\frak d_1$ and
$\frak d_2$. Contracting $\frak d_1$ instead corresponds to making an
elementary
modification of $\Bbb F_N$ and thus replacing it by $\Bbb F_{N\pm 1}$. As we
shall see, the symmetry between $\frak d_1$ and
$\frak d_2$ corresponds to the choice of either $\Delta$ or $\Delta - \kappa$.

The branch divisor $B'$ on the blowup $\Bbb F$ of
$\Bbb F_N$ is of the form $B+ \frak d_1+\frak d_2$, where $B$ is  a smooth
4-section which does not meet $\frak d_1$ or $\frak d_2$ and hence  $B\cdot
\frak e = 2$. Thus if we use the basis $\{\sigma, f, \frak d_2, \frak e\}$ for
$\operatorname{Pic}(\Bbb F)$, where $\sigma $ is the negative section of  $\Bbb
F_N$ and $f$ is the class of a fiber, viewed as curves on $\Bbb F$, it is easy
to see that  $$B = 4\sigma + (2k+1)f -4\frak e -2\frak d_2$$  for some odd
integer $k$ and that $$B+\frak d_1+\frak d_2 = 4\sigma + (2k+2)f -6\frak e -2
\frak d_2,$$ which is indeed divisible by 2. Note that we cannot say  {\it a
priori\/} that $B$ is irreducible. However it cannot be a union of a 3-section
and a section, since there are no sections of  $J^{m_1m_2}(S)$. In particular
$B\cdot \sigma \geq 0$. Thus if $\frak d_1$ is the proper transform of the
fiber on $\Bbb F_N$, and we assume that the negative section
$\sigma$ does not pass through the point of the original fiber that was blown
up, so that $\sigma \cdot \frak d_2 = 0$, then  $2k+1 \geq 4N$, or equivalently
$$k\geq 2N.$$
The same conclusion holds by a similar argument if $\sigma$ does pass through
the point that is blown up.

It is now easy to reverse this procedure. Begin with $\Bbb F_N$ and blow up a
point in a fiber. Then blow up the point of intersection of the exceptional
curve with the proper transform of the fiber. The result is a a
nonmiminal ruled surface $\Bbb F$ with a reducible fiber of the ruling of the
form $\frak d_1+2\frak e + \frak d_2$, where $\frak d_1$ is the proper
transform
of the fiber, $\frak d_2$ is the proper transform of the first exceptional
curve, and $\frak e$ is the second exceptional curve. Choose a smooth element
$B$ in the linear system $|4\sigma +(2k+1)f - 4\frak e -2\frak d_2|$, if any
exist, where $\sigma$ is the negative section of $\Bbb F_N$ and $f$ is a fiber.
The double cover of $\Bbb F$ branched along $B+\frak d_1+\frak d_2$ is then an
elliptic surface with a multiple fiber of multiplicity 2, bisections
corresponding to the pullbacks of sections of $\Bbb F_N$, and an involution
$\iota$.

Let us calculate $p_g(J^{m_1m_2}(S)) = p_g(S)$ in terms of $\Bbb F$ and $B$.
The
canonical bundle of $\Bbb F_N$ is given by $K_{\Bbb F_N} = -2\sigma -(N+2)f$.
Thus, recalling that $\frak d_2$ is the proper transform of the first
exceptional curve and that $\frak e$ is the second, we have
$$K_{\Bbb F} =  -2\sigma -(N+2)f +\frak d_2 +2\frak e.$$
As for the branch locus $B+\frak d_1+\frak d_2$, we have
$$B+\frak d_1+\frak d_2 = 2(2\sigma +(k+1)f -3\frak e -\frak d_2).$$
By standard formulas for double covers,
$$H^0(J^{m_1m_2}(S); K_{J^{m_1m_2}(S)}) = H^0(\Bbb F; \scrO _{\Bbb F}(K_{\Bbb
F} +
2\sigma +(k+1)f -3\frak e -\frak d_2)).$$
Now, using the calculations above, we have
$$K_{\Bbb F} +
2\sigma +(k+1)f -3\frak e -\frak d_2 = (k-N-1)f - \frak e.$$
Recalling that $f$ is linearly equivalent to $\frak d_1+2\frak e + \frak d_2$,
it is clear that
$$h^0((k-N-1)f - \frak e) = \cases k-N-1, &\text{if $k-N\geq 2$}\\
0, &\text{otherwise.} \endcases$$
Now if $k-N \leq 1$, since $k\geq 2N$ we must have $N\leq 1$. If $N=1$, then
$k=2$ and so $k-N-1=0$ in this case as well. If $N=0$, then $k=0$. In this
case $B= 4\sigma +f-4\frak e -2\frak d_2= 4\sigma +\frak d_1-2\frak e - \frak
d_2$, and it is easy to see that the effective curve $\sigma - \frak d_2-
\frak e$, which is
the proper transform of the unique element of $|\sigma|$ passing through the
point of the fiber which is blown up, satisfies
$$ (\sigma - \frak d_2-\frak e)\cdot (4\sigma +f-4\frak e -2\frak d_2) =-1.$$
Thus $|B|$ has the fixed component $\sigma - \frak d_2-\frak e$, which is a
section, and this case does not arise. So in all cases we have $h^0((k-N-1)f -
\frak e) = k-N-1$.

Thus we may summarize this discussion as follows:

\lemma{1.2} With notation and conventions as above,
$$p_g(S) =p_g(J^{m_1m_2}(S)) =  k-N-1. \qed$$
\endstatement

In particular, suppose that $p_g(S) = 0$. Since $k\geq 2N$, and the case $k=N =
0$ has been ruled out above, the only possibilities are $N = 0$, $k=1$ or
$N=1$, $k=2$. Thus $\Bbb F$ is the blowup of either $\Bbb F_0$ or $\Bbb F_1$,
and of course the two cases are elementary transformations of each other. Thus
we may assume that $k=2$ and $N=1$ in this case. Moreover the negative section
of $\Bbb F_1$ does not pass through the exceptional point in the blowup.

\section{2. Classification of stable bundles.}

We let $S$ be a nodal elliptic over $\Pee ^1$ with  exactly two multiple fibers
of multiplicities $2m_1$ and $m_2$ and let $\Delta$ be a divisor on $S$ with
$\Delta \cdot f = 2m_1m_2$. Fix an integer $c$. In this section we shall study
rank two vector bundles $V$ with $c_1(V) = \Delta$ and $c_2(V) = c$. We shall
also let $w = \Delta \bmod 2$ and $p = \Delta ^2 -4c$.

First  recall the following standard definition from [3]:

\definition{Definition 2.1} An ample line bundle $L$ on $S$ is {\sl
$(\Delta, c)$-suitable\/} or {\sl
$(w,p)$-suitable\/} if for all divisors $D$ on $S$ such that $-D^2
+ D\cdot \Delta \leq c$, either $f\cdot(2D - \Delta) = 0$ or
$$\operatorname{sign} f\cdot (2D - \Delta) =
\operatorname{sign}L\cdot (2D - \Delta).$$
\enddefinition

The following is  Lemma 3.3 of [3]:

\lemma{2.2} For all pairs $(\Delta, c)$, $(\Delta, c)$-suitable ample line
bundles exist. \qed
\endstatement

With this said, here is the rough classification of rank two
vector bundles $V$ with $c_1(V) = \Delta$ and $c_2(V) = c$ which
are stable with respect to a $c$-suitable line bundle $L$.

\theorem{2.3} Let $L$ be $(\Delta, c)$-suitable, and let $V$ be an $L$-stable
rank two vector bundle with  $c_1(V) = \Delta$ and $c_2(V) = c$. Then there
exist:  \roster
\item"{(i)}" A smooth irreducible curve $C$ and a birational map $C \to
\overline C \subseteq J^{m_1m_2}(S)$, where $\overline C$ is a bisection of
$J^{m_1m_2}(S)$ invariant under the involution $\iota$;
\item"{(ii)}" A divisor $D$ on $ T$, the minimal desingularization of the
normalization of $C\times _{\Pee^1}S$, such that $D\cdot f = m_1m_2$, where $f$
is a general fiber of $T\to C$, and moreover $D$ has the same restriction
to the generic fiber of $ T$ as the divisor induced by the section of
$J^{m_1m_2}(T)$ corresponding to the map $C \to J^{m_1m_2}(S)$;
\item"{(iii)}" A codimension two local complete intersection $Z$ and an exact
sequence
$$0 \to \scrO_{T}(D) \to \nu ^*V \to \scrO_{T}(\nu ^*\Delta -D)\otimes I_Z \to
0,$$
where $\nu \:T \to S$ is the natural degree two map.
\endroster
Moreover the bisection $\overline C$ and the double cover $T$ are
uniquely determined by the bundle $V$, and $D$ is determined by the bundle $V$
and the choice of a map $C \to J^{m_1m_2}(S)$. Finally, every rank two vector
bundle $V$ with $c_1(V) = \Delta$ and $c_2(V) = c$ satisfying
\rom{(i)}--\rom{(iii)} above is stable with respect to every $(\Delta,
c)$-suitable ample line bundle $L$. \endstatement
\proof First suppose that $V$ is $L$-stable. It follows from Theorem 4.3 of
[3] that the restriction of $V$ to the geometric fiber of $\pi$ is
semistable. More precisely,
let $\eta = \operatorname{Spec} k(\Pee^1)$ and let $\bar \eta =
\operatorname{Spec} \overline{k(\Pee^1)}$, where $\overline{k(\Pee^1)}$ denotes
the algebraic closure of $k(\Pee^1)$. Let $V_{\bar \eta}$ denote the pullback
of
$V$ to the curve $S_{\bar \eta}$ which is the geometric fiber of $\pi$. Then
$V_{\bar \eta}$  is semistable.  By the classification of rank two bundles
on an elliptic curve, $V_{\bar \eta}= L_1 \oplus L_2$, where each $L_i$ is a
line bundle over $S_{\bar \eta}$ of degree $m_1m_2$ and $L_1\otimes L_2$
corresponds to the restriction of $\Delta$ to $S_{\bar \eta}$. The Galois group
$\operatorname{Gal}(\overline{k(\Pee^1)}/k(\Pee ^1)$ permutes the set $\{L_1,
L_2\}$. This action cannot be trivial, since otherwise $L_i$ would be rational
over $k(\Pee ^1)$ and then $S$ would have an $m_1m_2$-section. Thus the fixed
field of the subgroup of  $\operatorname{Gal}(\overline{k(\Pee^1)}/k(\Pee ^1)$
which operates trivially on $\{L_1,
L_2\}$ defines a degree two extension of $k(\Pee^1)$, corresponding to a
morphism $C\to \Pee ^1$. Setting $T $ to be the minimal resolution of the
normalization of $C\times _{\Pee ^1}S$, there is a section of $J^{m_1m_2}(T
)$ defined by $L_1$, say. The image of this section in $J^{m_1m_2}(S)$ is then
the
bisection $\overline C$. By construction $\overline C$ is invariant under the
involution $\iota$. Let $\nu \:T \to S$ be the natural degree two map.

The inclusion $L_1 \to V_{\bar \eta}$ induces a sub-line bundle $\scrO_{T}(D)
\to \nu ^*V$, which we may assume to have torsion free cokernel. Since $L_2
\neq
L_1$, it is clear that this sub-line bundle is unique. The quotient is then
necessarily of the form $\scrO_{ T}(\nu ^*\Delta -D)\otimes I_Z$.

It remains to prove that every $V$ satisfying the above description is indeed
$L$-stable. It follows from  (iii) that $ V_{\bar
\eta}$ is an extension of two line bundles of degree $m_1m_2$ and is therefore
semistable. Again using Theorem 4.3 of [3], $V$ is $L$-stable.
\endproof

Next we discuss the meaning of the scheme $Z$ and the bisection $\overline C$.
The following is the analogue of Lemma 1.11 in [4] and is proved in exactly
the same way:

\lemma{2.4} Let $f$ be a smooth fiber of $\pi$ and let $g$ be a component of
$\nu^{-1}(f)$. Then $\operatorname{Supp}Z\cap g \neq \emptyset$ if and
only if $V|f$ is unstable. In particular, if $\nu$ is not branched
over $f$, so that $\nu^{-1}(f) = g\cup g'$, then
$\operatorname{Supp}Z\cap g \neq \emptyset$ if and only if
$\operatorname{Supp}Z\cap g' \neq \emptyset$.
\qed
\endstatement

Next we turn to the section $\overline C$. Since $\overline C$ is invariant
under $\iota$, its proper transform on the blowup of $J^{m_1m_2}(S)$ at the two
isolated fixed points of $\iota$ is the pullback of a section $A'$ of $\Bbb F$,
which in turn induces a section $A$ of $\Bbb F_N$. We shall use throughout
the notation and conventions of the previous section. Notice that the section
$A'$ meets the reducible fiber either along $\frak d_1$ or $\frak d_2$,  the
two
components of multiplicity one. Here $A'\cdot \frak d_1 =1$ and $A'\cdot
\frak d_2=0$  if $A$ does not pass through the point of the corresponding fiber
of $\Bbb F_N$ which was blown up, and $A' \cdot \frak d_2 =1$ and
$A'\cdot \frak d_1=0$ in
the remaining case. Since the branch locus of the map $J^{m_1m_2}(S) \to \Bbb
F$
consists of $B+\frak d_1+\frak d_2$, we see that $A'$ always passes through the
branch locus over the point corresponding to $F_2$. Of course, this is also
clear from the picture on $J^{m_1m_2}(S)$: since $\overline C$ is a bisection,
$\overline C \cdot f =2$ and therefore $\overline C \cdot F_2=1$. It follows
that $\overline C$ is smooth at the point of intersection with $F_2$, that the
intersection is transverse, and that the natural map $\overline C \to \Pee ^1$
is always branched at the point corresponding to $F_2$. A similar
statement will hold for the map $C\to \Pee ^1$. This fact is the fundamental
difference between the case studied in this paper and the case of trivial
determinant studied in [2] and [4].

Since $A'$ always meets $\frak d_1$ or $\frak d_2$, it follows that the inverse
image of $A'$ in the blowup of $J^{m_1m_2}(S)$ always meets exactly one of the
two
exceptional curves, and in fact meets it transversally at one point. Thus the
inverse image of $A'$ is the proper transform of $\overline C$, and therefore
$$(\overline C)^2 = 2(A')^2 + 1.$$

Let us consider the section $A$ of $\Bbb F_N$ in more detail. Either $A=
\sigma$ or $A\in |\sigma + (N+s)f|$ for a uniquely specified nonnegative
integer $f$. Moreover either $A$ does not pass through the point on $\Bbb F_N$
which is the image of the exceptional divisor, in which case $A'\cdot \frak
d_2 =0$, or it does, in which case $A' \cdot \frak d_2=1$. The following lemma
relates the odd integer $\Delta ^2 - 4c=p$ to the invariants of $V$:

\lemma{2.5} With notation as above, denote by the exceptional point the point
of $\Bbb F_N$ which is blown up under the morphism $\Bbb F \to \Bbb F_N$. Then,
if we set $p = p_1(\ad V) = \Delta ^2 - 4c$,
$$-p = \cases 4k-6N+1+2\ell(Z) +\delta, &\text{if $A=\sigma$ does not pass
through}\\
{} &\text{the exceptional point;}\\
4s+4k-2N+1+2\ell(Z) +\delta, &\text{if $A\in |\sigma +(N+s)f|$}\text{does not
pass} \\   {} &\text{through the exceptional point;}\\
4k-6N-1+2\ell(Z) +\delta, &\text{if $A=\sigma$ passes through}\\
{} &\text{the
exceptional point;}\\
4s+4k-2N-1+2\ell(Z) +\delta, &\text{if $A\in |\sigma +(N+s)f|$ passes
through}\\
{}&\text{the
exceptional point.}
\endcases$$
Here $\delta$ is a nonnegative integer which is zero
if the map $C\to \Pee ^1$ is not branched over any point
corresponding to a singular nonmultiple fiber of $\pi \: S \to \Pee ^1$.
\endstatement
\proof Since $c_1^2(\nu ^*V) = 2c_1^2(V)=2\Delta ^2$ and $c_2(\nu
^*V) = 2c_2(V) = 2c$, it will suffice to work with $\nu ^*V$. Clearly
$$c_2(\nu ^*V) = -D^2 + D\cdot \nu ^*\Delta + \ell (Z).$$
Thus
$$2(4c - \Delta ^2) = -(2D- \nu ^*\Delta )^2 + 4\ell(Z).$$
Now we can write
$$2D- \nu ^*\Delta = D-(\nu ^*\Delta -D),$$
where both $D$ and $\nu ^*\Delta -D$ naturally correspond to sections of
$J^{m_1m_2}(T)$. In fact, if $\varphi\: J^{m_1m_2}(T) \to J^{m_1m_2}(S)$ is
the obvious map, then the bisection $\overline C$ satisfies $\varphi
^*\overline C = C_1 + C_2$, where $C_1$ and $C_2$ are sections of
$J^{m_1m_2}(T)$ corresponding to the divisors $D$ and $\nu ^*\Delta
-D$ on $T$. An argument essentially identical to the proofs of Claim
1.17 and 1.18 in Chapter 7 of [4] shows that there is a nonnegative integer
$\delta$ such that
$$ -(2D- \nu ^*\Delta )^2  = -(C_1-C_2)^2 + 2\delta.$$
Moreover $\delta = 0$ if the map $C\to \Pee ^1$ is not branched over any point
corresponding to a singular nonmultiple fiber of $\pi \: S \to \Pee ^1$. Thus
we must calculate $(C_1-C_2)^2$. But, using the fact that $C_i$ is a section of
$J^{m_1m_2}(T)$, we have
$$(C_i)^2 = -2(1+p_g(S)) = -2(k-N).$$
Moreover $C_1+C_2 = \varphi ^*\overline C$. Thus
$$\align
(C_1-C_2)^2 &= 2(C_1)^2 + 2(C_2)^2 - (C_1+C_2)^2\\
&=-8(k-N) -2(\overline C)^2\\
&=-8(k-N) -4(A')^2 - 2.
\endalign$$
Clearly we have
$$(A')^2 = \cases -N, &\text{if $A=\sigma$ does not pass through the
exceptional point;}\\
N+2s, &\text{if $A\in |\sigma +(N+s)f|$ does not pass through}\\
{} &\text{the
exceptional point;}\\
-N-1, &\text{if $A=\sigma$ passes through the
exceptional point;}\\
N+2s-1, &\text{if $A\in |\sigma +(N+s)f|$ passes through}\\
{}&\text{the
exceptional point.}
\endcases$$

Putting these formulas together gives the statement of the lemma.
\endproof

Using Lemma 2.5, the inequality $k\geq 2N$ and the fact that if $N=0$ then
$k\geq 1$, whereas if $N=1$ and $k=2$ then the section $\sigma$ does not pass
through the exceptional point, we can easily deduce the following slight
strengthening of Bogomolov's inequality in our case:

\corollary{2.6} We have the following inequality for $-p$:
$$-p\geq \cases 4p_g(S) -2N +5 &\text{if $A=\sigma$ does not pass through the
exceptional point;}\\
4p_g(S) +2N +5 &\text{if $A\in |\sigma +(N+s)f|$ does not pass through}\\
{} &\text{the
exceptional point;}\\
 4p_g(S) -2N +3 &\text{if $A=\sigma$ passes through the
exceptional point;}\\
4p_g(S) +2N +3  &\text{if $A\in |\sigma +(N+s)f|$ passes through}\\
{}&\text{the
exceptional point.}
\endcases$$
In all cases $-p \geq 2p_g(S) +1$, and if $p_g(S)=0$, then $-p\geq 3$. \qed
\endstatement

\section{3. A Zariski open subset of the moduli space.}

Our goal in this section is to prove the following theorem:

\theorem{3.1} Let $V$ be an $L$-stable rank two bundle on $S$. Suppose that
\roster
\item"{(i)}" The associated bisection $\overline C$ of $J^{m_1m_2}(S)$ is
smooth,
or equivalently that $\overline C = C$, and the image of $\overline C$ in $\Bbb
F$ is not the proper transform of $\sigma$; \item"{(ii)}" The map $C \to \Pee
^1$ is not branched at any point corresponding to a singular fiber of $\pi$ or
at the multiple fiber of odd multiplicity $m_2$;
\item"{(iii)}" The scheme $Z$ on the associated double cover $T$ is empty, and
thus there is an exact sequence
$$0 \to \scrO _T(D) \to \nu ^*V \to \scrO _T(\nu ^*\Delta - D) \to 0.$$
\endroster
Then $V = \nu _*\scrO _T(D+F)=\nu _*\scrO _T(\nu ^*\Delta -D)$. In particular
$V$ is uniquely determined by the associated section $A$ of $\Bbb F_N$ and the
divisor $D$ on $T$. Finally $V$ is a smooth point of its moduli space, which is
of dimension $-p -3\chi (\scrO_S)$ at $V$. \endstatement
\medskip

It is clear that the conditions above are equivalent to assuming that $A'$
meets the branch locus $B$ transversally, and that no point of intersection
lies over a point of $\Pee ^1$ corresponding to a singular nonmultiple fiber
or to the multiple fiber of multiplicity $m_2$, and that $Z=\emptyset$.

The proof of (3.1) will proceed along lines very similar to the proof of
Theorem 1.12 in Chapter 7 of [4], and we shall simply sketch some of the
details.

Let $A$ be the section of $\Bbb F_N$ corresponding to $A'$. By assumption $A
\neq \sigma$. Let $r$ be the nonnegative integer such that $A\in |\sigma
+(N+r)f|$. If the section $A$ does not pass through the exceptional point of
the
blowup, then  $$(A')\cdot (B+\frak d_1+\frak d_2) = (\sigma +(N+r)f)\cdot
(4\sigma +(2k+2)f -6\frak e -2\frak d_2)= 4r+2k+2.$$ Of these points, one
corresponds to the intersection $A'\cdot \frak d_1$, and so the branch divisor
of the map $T \to S$ is $(4r+2k+1)f$, where $f$ is a general fiber of $\pi$.
This
divisor is even since $f$ is divisible by $2$, and we set $G =(4r+2k+1)f/2$.
Likewise, if $A$ does pass through the exceptional point of the blowup, then
$$(A')\cdot (B+\frak d_1+\frak d_2) = (\sigma +(N+r)f-\frak d_2-\frak e)\cdot
(4\sigma +(2k+2)f -6\frak e -2\frak d_2)= 4r+2k.$$ In this case we set $G=
(4r+2k-1)f/2$. Let $F$ be the branch divisor in $T$, so that $\nu ^*G\equiv F$.
Thus $F = (4r+2k+1)f$ or $(4r+2k-1)f$. For future reference, let us also record
the genus of $C$:

\lemma{3.2} Let $C$ satisfy \rom{(i)} and \rom{(ii)} of \rom{(3.1)}. Then
$$g(C) = \cases 2r+k, &\text{if $A$ does not pass through the exceptional
point;}\\
2r+k-1, &\text{if $A$ passes through the exceptional point.}
\endcases$$
\endstatement
\proof The map $C\to \Pee ^1$ is branched at $A'\cdot (B+\frak d_1+\frak d_2) =
4r+2k+2$ points if $A$ does not pass through the exceptional point, and $4r +
2k$ points otherwise. The lemma now follows from the Riemann-Hurwitz formula.
\endproof

Now $\det \nu
_*\scrO_T(D) = \nu _*D - G$.
Clearly $\nu _*D$ and $\Delta$ have the same restriction to the generic fiber.
Arguing as in Chapter 7, (1.20) of [4], there is an
injective map  $$\nu _* \scrO _T(D+F) \to V.$$
Set $W = \nu _* \scrO _T(D+F)$. Our goal will be to show that $W = V$.
Preliminary to this goal we shall analyze $W$ and the map $W\to V$.
As a divisor
class $\det W = \nu _*D-G +2G= \nu _* D + G$. In addition there is an effective
divisor $E$ such that $(\det W)^{-1} \otimes \det V = \scrO _S(E)$. Thus $\det
W
= \Delta - E$ and $E$ has trivial restriction to the generic fiber, so that $E$
is a union of fibers (possibly including the reductions of the multiple
fibers).  Moreover $$\nu _*D = \Delta - G - E.$$
Set $E' = \nu ^*E$. We have
$$D + \iota ^*D = \nu ^*\nu _*D = \nu ^*\Delta -F -E',$$
and therefore
$$\iota ^*D = \nu ^*\Delta - D -F -E'.$$
We can thus write
$$W = \nu _*\scrO _T(D+F) = \nu _* \scrO _T(\nu ^*\Delta -D
-E').$$
Using the fact that there is a surjection from $\nu ^*W$ to $\scrO _T(\nu ^*
\Delta -D -E')$, we conclude that there is an exact sequence
$$0 \to \scrO _T(D) \to \nu ^*W \to \scrO _T(\nu ^*
\Delta -D -E') \to 0.$$

Comparing this sequence to the defining exact sequence for $\nu ^*V$ and
arguing as in (1.24) of Chapter 7 of [4], we may conclude:

\lemma{3.3} Let $Q =V/W$. Then $\nu ^*Q \cong [\scrO _T/\scrO_T(-E')] \otimes
\scrO_T(\nu ^*\Delta -D)$.\qed
\endstatement
\medskip

Our goal now is to prove the following:

\lemma{3.4} In the above notation, $E'=0$. Thus $Q=0$ and $V = W =\nu _*
\scrO _T(D+F)$ where $\nu _* D = \Delta -G$.
\endstatement
\medskip

We begin with the following construction.
Let $e$ be a component of the support of $E'$, and write $E' = ae + E''$, where
$E''$ is effective and disjoint from $e$ and $a>0$. If $e$ is not the multiple
fiber of multiplicity $m_1$ on $T$, then either $\nu$ is unbranched over $e$ or
$e$ is a smooth fiber. In either of these cases $\nu$ induces an isomorphism
from
$e$ to $\nu (e)=f$, and we shall identify $\nu (e)$ with $e$. In the remaining
case $e=F_{m_1}$ is the multiple fiber of multiplicity $m_1$ and $\nu$ is an
\'etale double cover.  There is then the following analogue of (1.25) of
Chapter
7 of [4]:

\lemma{3.5} There is a subsheaf $Q_0$ of $\nu _*Q$ which is isomorphic to
\roster
\item"{(i)}" $\scrO _e(-(a-1)e+\nu ^*\Delta -D),$
viewed as a sheaf on $\nu (e) = f$, if $e\neq F_{m_1}$;
\item"{(ii)}" A line bundle on $F_{2m_1}$ such that
$$\nu ^*Q_0 \cong
\scrO _{F_{m_1}}(-(a-1)F_{m_1}+\nu ^*\Delta -D),$$
in case $e=F_{m_1}$.
\endroster
\endstatement
\proof The argument in case $e\neq F_{m_1}$ runs as in (1.25) of Chapter 7 of
[4]. If $e = F_{m_1}$, then $\nu ^*Q$ contains the subsheaf
$$Q_0' = [\scrO _T(-(a-1)e/\scrO_T(-ae)] \otimes
\scrO_T(\nu ^*\Delta -D),$$
which is a line bundle on $F_{m_1}$.
The vector bundle $\nu _*Q_0'$ is a rank two vector bundle on $F_{2m_1}$ with
$\deg (\nu _*Q_0') = m_1m_2$. Consider the rank two vector bundle $\nu ^* \nu
_*Q_0'$ on $F_{m_1}$. Its determinant is $Q_0'\otimes \iota^*Q_0' =
(Q_0')^{\otimes 2}$ (recall that $\nu ^*\Delta -D$ is fixed under the
involution) and there is a surjective map $\nu ^* \nu _*Q_0' \to Q_0'$. Thus
there is an exact sequence $$0 \to Q_0' \to \nu ^* \nu _*Q_0' \to Q_0' \to 0.$$
It follows that $\nu ^* \nu _*Q_0'$ is semistable and is either $Q_0'\oplus
Q_0'$ or the nontrivial extension of $Q_0'$ by $Q_0'$. On the other hand,
$\nu _*Q_0'$ is either a direct sum of line bundles, say $Q_0 \oplus Q_1$
for two line bundles $Q_i$ or a nontrivial extension of a line bundle
$Q_0$ by $Q_0$. In the first case we must have $\nu ^*Q_i = Q_0'$ and in the
second $\nu ^*Q_0 = Q_0'$. In either case there is a subbundle $Q_0$ of $\nu
_*Q$
as desired.
\endproof

To prove Lemma 3.4, we shall assume that $E'\neq 0$ and derive a contradiction.
Again, the argument will be very similar to the argument given in [4]. It will
suffice to show that $\dim \operatorname{Ext}^1(Q_0, W) \leq 1$. We also have
that
$$\operatorname{Ext}^1(Q_0, W) = H^0(W\otimes \scrO _S(e)\otimes Q_0^{-1}).$$
The case where $e$ does not lie over the branch locus of $C\to \Pee ^1$ follows
exactly as in [4]. The case where $e$ is a smooth fiber in the branch locus
also follows by these methods provided we can show that $R^1\rho _*\scrO
_T(2D-\nu ^*\Delta )$ has length one at the point of $C$ corresponding to $e$.
This is a local calculation, which we shall leave to the reader; it uses the
fact that $A'$ meets the branch locus transversally and can in fact be deduced
from the global argument in [4], proof of Lemma 1.19 of Chapter 7.

There remains the new case, where $e = F_{m_1}$. In this case, the natural map
$$W = \nu _* \scrO _T(\nu ^*\Delta -D -E')\to \nu _* \scrO _{F_{m_1}}(\nu
^*\Delta -D
-aF_{m_1})$$
is surjective, as one can see from applying the surjective map $\nu _*$ to the
exact sequence
$$0 \to \scrO _T (\nu ^*\Delta -D -E'-F_{m_1})\to
\scrO _T(\nu ^*\Delta -D -E')\to \scrO _{F_{m_1}}(\nu ^*\Delta -D
-aF_{m_1})\to 0.$$
It follows that $W|F_{2m_1} =  \nu _* \scrO _{F_{m_1}}(\nu
^*\Delta -D -aF_{m_1})$. Thus
$$\gather
H^0(W\otimes \scrO _S(F_{2m_1})\otimes Q_0^{-1}) =
H^0(\nu _*\scrO _{F_{m_1}}(\nu ^*\Delta -D
-aF_{m_1}) \otimes \scrO _S(F_{2m_1})\otimes Q_0^{-1}) \\
= H^0(\nu _*\big[\scrO _{F_{m_1}}(\nu ^*\Delta -D
-aF_{m_1}) \otimes \nu ^*\scrO _S(F_{2m_1})\otimes \nu ^*Q_0^{-1}\bigr]) =
H^0(\scrO _{F_{m_1}}), \endgather$$
where we have used the fact that $\nu ^*Q_0 = \scrO _{F_{m_1}}(\nu ^*\Delta -D
-(a-1)F_{m_1})$.
Hence
$$\dim\operatorname{Ext}^1(Q_0, W) = h^0(W\otimes \scrO _S(F_{2m_1})\otimes
Q_0^{-1}) = 1$$
as desired.
\qed
\medskip

We see that we have proved all of Theorem 3.1 except the statement about the
smoothness of the moduli space, which follows from:

\lemma{3.6} Suppose that $V= \nu _*\scrO _T(D+F)$ as above. Then
$V$ is good, in the terminology of \rom{[2]}. In other words, $H^2(S;
\operatorname{ad} V) = 0$.
\endstatement
\proof It suffices to show that $\dim \operatorname{Hom}(V, V\otimes K_S) =
h^0(K_S)$. Now
$$\align
\operatorname{Hom}(V, V\otimes K_S) &= \operatorname{Hom}(V, \nu _*\Big(\scrO
_T(\nu ^*\Delta -D)\otimes \nu ^*K_S\Bigr)) \\
&= \operatorname{Hom}(\nu ^*V, \scrO
_T(\nu ^*\Delta -D)\otimes \nu ^*K_S).
\endalign$$
Using the defining exact sequence for $\nu ^*$, there is an exact sequence
$$0 \to H^0(\nu ^*K_S) \to \operatorname{Hom}(\nu ^*V, \scrO
_T(\nu ^*\Delta -D)\otimes \nu ^*K_S) \to H^0(\scrO _T(\nu ^*\Delta -2D)\otimes
\nu ^*K_S).$$
Since $K_S$ is a rational multiple of the fiber and $\nu ^*\Delta -2D$ is
nontrivial on the generic fiber, the term $H^0(\scrO _T(\nu ^*\Delta
-2D)\otimes
\nu ^*K_S)$ is zero. Thus
$\dim \operatorname{Hom}(V, V\otimes K_S) = h^0(\nu ^*K_S)$. Using the
isomorphism $H^0(\nu ^*K_S) \cong H^0(K_S) \oplus H^0(K_S(-G))$, it suffices to
show that $ H^0(K_S(-G)) = 0$. Now $K_S= \scrO_S((k-N-2)f + (2m_1-1)F_{2m_1} +
(m_2-1)F_{m_2})$. Also $G = (2r+k \pm 1/2)f$. Thus it suffices to observe that
the linear system
$$|(-N -2r -2\pm 1/2)f + (2m_1-1)F_{2m_1} +
(m_2-1)F_{m_2})|$$
is empty.
\endproof

Next let us describe the subset of the moduli space consisting of bundles $V$
which satisfy the hypotheses of Theorem 3.1. We begin by reversing the
procedure outlined above. Fix the section $A'$, which is generic in the sense
of Theorem 3.1: it meets $B$ transversally and no point of intersection
corresponds to a singular nonmultiple fiber or to the multiple fiber of odd
multiplicity. The section $A'$ determines the bisection $C = \overline C$, and
thus a double cover $\nu \: T \to S$ together with an elliptic fibration $\rho
\:
T\to C$ and a divisor $D_0$, well-defined on the generic fiber.
Moreover by construction $\nu _*D_0$ and $\Delta$ have the same restriction on
the generic fiber, and thus differ by a multiple of $\kappa$. It is easy to see
that changing $D_0$ by a sum of fiber components on $T$ replaces $\nu _*D_0$ by
an arbitrary even multiple of $\kappa$. Thus we may assume that we have
$\nu _* D_0 = \Delta - G$ or $\nu _* D_0 = \Delta - \kappa - G$. It is an
exercise in the formulas of the preceding section to see that we have
$$\Delta ^2 \equiv 2(A')^2 + 1 \mod 4.$$
Additionally
$$2(A')^2 + 1 \equiv \cases 2N+1 \mod 4, &\text{if $A$ does not pass through
the
exceptional point;}\\
2N-1 \mod 4 , &\text{otherwise.}
\endcases$$
 Thus the symmetry between the possibility that $A$ does or does not pass
through the exceptional point, which is essentially the choice of blowing $\Bbb
F$ down to $\Bbb F_N$ or $\Bbb F_{N\pm 1}$, reflects the choice of $\Delta$ or
$\Delta - \kappa$, which in turn reflects $\Delta ^2\bmod 2$, or equivalently
the dimension of the moduli space mod 2.

Having made one choice for a  line bundle $D_0$ on the double cover $T$, where
$D_0$ is specified on the generic fiber of $T$ and satisfies $\nu _*D_0 =
\Delta
- G$, or $\nu _*D_0 = \Delta - \kappa -G$, the  possibilities for $D$
are given by the next lemma.

\lemma{3.7} Given the double cover $T\to S$, the set of all $D$ whose
restriction
to the generic fiber equals $D_0$ and which satisfy $\nu _*D = \Delta
- G$, or $\nu _*D = \Delta - \kappa -G$ is a principal homogeneous space over
$\operatorname{Pic}^\tau T$, which in turn is an extension of the Jacobian
$J(C)$ by a cyclic group of order $m_2$. Moreover this principal
homogeneous space is nonempty for exactly one of the two choices for
$\nu _*D$ above.
\endstatement
\proof By the remarks preceding the lemma, there exists a $D_0$ with
$\nu _*D_0 = \Delta
- G$, or $\nu _*D_0 = \Delta - \kappa -G$, and only one of these possibilities
can hold. If $D$ has the same restriction to the generic fiber as $D_0$ and
$\nu _*D = \nu _*D_0$, then $D-D_0$ has trivial restriction to the generic
fiber and $\nu _*(D-D_0) = 0$. The first condition says that $D-D_0$ is of the
form $\rho ^*\lambda \otimes \scrO_T(aF_{m_1} + bF'_{m_2} + cF''_{m_2})$, where
$F_{m_1}$ is
the multiple fiber of multiplicity $m_1$ lying above $F_{2m_1}$ and $F'_{m_2}$,
$F_{m_2}''$ are the two multiple fibers of multiplicity $m_2$ lying over
$F_{m_2}$. We may further assume that $0\leq a <m_1$ and that $0\leq b < m_2$,
$0\leq c <m_2$.  Here $\lambda$ is a line bundle of degree $d$ on $C$. Thus
$$\nu _* (D-D_0)  = df + 2aF_{2m_1} + (b+c)F_{m_2}.$$
It is easy to see that this divisor is trivial if and only if $d=0$, $a=0$,
and $b\equiv -c\mod m_2$. Thus, there is a natural identification of the set of
all $D$ (given the fixed divisor $D_0$) with $J(C) \times \Zee/m_2\Zee$.
\endproof

We now assume that $A$ is not the negative section of $\Bbb F_N$ and write
$A\in |\sigma +(N+r)f|$ where $r\geq 0$. The dimension of the linear system
$A'$ is then equal to $n$, where   $$n = \cases N+2r +1, &\text{if $A$ does not
pass through the  exceptional point;}\\ N+ 2r, &\text{otherwise.}
\endcases$$
We also let $g = g(C)$ be the genus of the bisection $C$, as given in Lemma
3.2. Note that
$$n+g = \cases 4r +N + k + 1  &\text{if $A$ does not
pass through the  exceptional point;}\\
4r + N + k -1 &\text{otherwise.}
\endcases$$
In both cases, comparing this with the formula for $-p$ given in Lemma 2.5 and
using the fact that $1+ p_g(S) = k-N$, we see that
$$-p-3\chi (\scrO _S) = g+n.$$

Note finally that the moduli space will be nonempty provided that $-p \geq
4p_g(S) +2N+3$. Since $p_g(S) = k-N -1 \geq N-1$, the moduli space will be
nonempty as long as
$$-p \geq 6p_g(S) + 5.$$

Arguing as in Theorem 1.14 of Chapter 7 of [4], we obtain the following:

\theorem{3.8} Let $p$ be an odd negative integer, and choose $w\in H^2(S;
\Zee/2\Zee)$ such that $w = \Delta \bmod 2$ or $w= \Delta - \kappa \bmod 2$,
and that $w^2 \equiv p\mod 4$. Let  $L$ be a $(w,p)$-suitable ample
line bundle on $S$. Let $\frak M = \frak M(S,L; w,p)$ denote the moduli
space of $L$-stable rank two bundles on $S$ with $w_2(V) = w$ and $p_1(\ad V) =
p$. Then for all $p$ such that $-p \geq 6p_g(S)+5$, $\frak M$ contains a
nonempty
Zariski open subset $M$ corresponding to vector bundles $V$ satisfying the
hypotheses of Theorem \rom{3.1}. The set $M$ is smooth of dimension $-p-3\chi
(\scrO _S)$. Moreover, there is a holomorphic map from $M$ to a Zariski open
subset $U \subseteq \Pee ^n$ and the fibers are isomorphic to $m_2$ copies of a
complex torus of dimension $g$. \qed \endstatement

Let us finally consider the case where $p_g(S) =0$ and $-p=3$, the case of a
moduli space of expected dimension zero. In this case we fix
$N=1$ and $k=2$, and the negative section does not pass through the exceptional
point. The Chern class calculations of Lemma 2.5 show that we must have $\ell
(Z) = \delta =0$ and $A$ must be the negative section $\sigma$ of $\Bbb F_1$.
Note that $\sigma \cdot (B+\frak d_1+\frak d_2) = 2$, and the intersection must
be transverse since $\sigma$ cannot split into a union of two sections in the
double cover. Thus $C= \overline C \to \Pee ^1$ is branched at two points, so
that $C = \Pee ^1$ again. Assuming as we may that the multiple fiber of odd
multiplicity does not correspong to a branch point, we see that the arguments
of Theorem 3.1 go through to show that there are exactly $m_2$ vector bundles
$V$ whose associated section is $\sigma$. (Here, in case the intersection
point of $\sigma$ with $B$ corresponds to a singular nonmultiple fiber, we
must use the more detailed analysis of [2] (5.12) and (5.13) to see that the
section $\sigma$ and the line bundle on $T$ determine $V$.) Each of these is a
smooth point of the moduli space, by a slight modification of the proof of
Lemma
3.6 (in this case $G= f/2$ and $K_S = \scrO_S(-f + (2m_1-1)F_{2m_1} +
(m_2-1)F_{m_2})$). Summarizing, then:

\theorem{3.9} In case $p_g(S) = 0$, the moduli space corresponding to $-p = 3$
consists of $m_2$ reduced points. \qed
\endstatement

\section{4. Calculation of the leading coefficient.}

Fix $w$ and $p$ with $w^2 \equiv p \mod 4$, and let $\frak M = \frak M(S,L; w ,
p)$ denote the  moduli space of $L$-stable rank two bundles on $S$ with $w_2(V)
= w$ and $p_1(\ad V) = p$. Let $d$ be the (complex) dimension of $\frak M$:
$$d = \cases 4r +N + k + 1  &\text{if $A$ does not
pass through the  exceptional point}\\
4r + N + k -1 &\text{otherwise,}
\endcases$$
where the section $A\in |\sigma +(N+r)f|$.
With this notation we see that $2n = d-p_g(S)$ and that $g= d-n$.
Finally let us recall from (3.4) of [3] that in case $p_g(S)=0$ there is a
unique chamber of type $(w,p)$ which contains $\kappa$ in its closure, called
the {\sl $(w, p)$-suitable\/} chamber. We can now state the main result of
this paper:

\theorem{4.1} For $p_g(S) >0$, let $\gamma _{w, p}(S, \beta)$ be the Donaldson
polynomial corresponding to the $SO(3)$-bundle $P$ over $S$ with $w_2(P) =
\Delta \bmod 2$ and $p_1(P) = p$, and $\beta$ is a choice of orientation
agreeing with the usual complex orientation for $\frak M$. If $p_g(S)=0$, let
$\gamma _{w, p}(S, \beta)$  be the corresponding Donaldson polynomial for
metrics whose associated self-dual harmonic $2$-form lies in the
$(w,p)$-suitable chamber. Then, writing  $\gamma _{w, p}(S, \beta)$ as a
polynomial in $\kappa _S$ and $q_S$, say  $\gamma _{w, p}(S, \beta)= \sum
_{i=0}^{[d/2]}a_iq_S^i\kappa _S^{d-2i}$, we have, for all $p$ with $-p\geq
2(4p_g+2)$, $a_i= 0$ for $i>n$ and  $$a_n =
\frac{d!}{2^nn!}(2m_1m_2)^{p_g(S)}m_2.$$ \endstatement \medskip

We first remark that the assumption that $-p\geq
2(4p_g+2)$ implies that $\frak M$ is nonempty and contains a smooth Zariski
open subset as described in Theorem 3.8. Indeed $-p$ is an odd integer greater
than $8p_g+4$, and so $-p \geq 8p_g+5 \geq 6p_g+5$. Thus we are in the range of
Theorem 3.8.

Let $X= X_{w,p}$ denote the Uhlenbeck compactification associated to $\frak M$
[1], [4]. The orientation of $\frak M$ induces a fundamental class of $X$.
There is a $\mu$-map $H_2(S) \to H^2(\frak M)$, which roughly speaking is given
by taking slant product with the class $-p_1(P)/4$, where $P$ is the universal
$SO(3)$-bundle over $S\times \frak M$. If $P$ lifts to a holomorphic bundle
$\Cal V$  over $S\times \frak M$, then $p_1(P) = c_1^2(\Cal V) - 4c_2(\Cal
V)=p_1(\ad \Cal V)$. The classes $\mu (\alpha) \in H^2(\frak M)$ extend
uniquely
to classes in $H^2(X)$, which we shall also denote by $\mu (\alpha)$. The
Donaldson polynomial is then defined by taking cup products of the $\mu
(\alpha)$ and evaluating on the fundamental cycle of $X$.

Arguing as in [4], it is a combinatorial exercise to deduce Theorem 4.1 from
the
following:

\theorem{4.2} Let $\mu$ denote the $\mu$-map associated to the Uhlenbeck
compactification $X$ of $\frak M$. Then, using the complex orientation to
identify $H^{2d}($X$; \Zee)\cong \Zee$, we have,
for all $p$ with $-p\geq 2(4p_g+2)$ and for all $\Sigma \in
H_2(S)$,   $$\mu (f)^m\cup \mu
(\Sigma)^{d-m} = \cases 0, &\text{if $m>n$}\\ (d-n)!(2m_1m_2)^{d-n}m_2(\Sigma
\cdot \kappa _S)^{d-n}, &\text{if $m=n$.} \endcases$$
\endstatement

Here the factor $(2m_1m_2)^{p_g(S)}$ appearing in Theorem 4.1 arises from the
fact
that $\mu (f)^n = (2m_1m_2)^n\mu (\kappa)$ and that $d-2n = p_g(S)$.

In order to prove Theorem 4.2, we shall
introduce geometric divisors which will represent the cohomology class
$\mu(f)$. Fix a smooth fiber $f$. Then there are exactly four line bundles
$\theta$ of degree $m_1m_2$ on $f$ such that $\theta ^{\otimes 2} = \scrO
_f(\Delta)$. Each line bundle $\theta$ corresponds to a point of intersection
of $f$ with the branch divisor $B\subset \Bbb F$. If $V$ is a rank two vector
bundle over $S$ with $c_1(V) = \Delta$, then   $(V|f)\otimes (\theta)^{-1}$ is
a
vector bundle over $f$ with trivial determinant. Thus by the Riemann-Roch
theorem $\chi (f; (V|f)\otimes (\theta)^{-1}) =0$.

For each integer $c$, fix a $(\Delta, c)$-suitable line bundle $L$. Given an
integer $b\leq c$, we define $\frak M_b$ to be the
moduli space of $L$-stable rank two bundles $V$ with $c_1(V) = \Delta$ and
$c_2(V) = b$. Given $f$ and $\theta$ as above, define the divisor
$\Cal Z_b(f, \theta)$ as a set by:
$$\Cal Z_b(f, \theta) = \{\, V \in \frak M_b: h^0(f; (V|f)\otimes
(\theta)^{-1}) \neq 0\,\}.$$
A calculation with the Grothendieck-Riemann-Roch theorem  as in the
proof of Proposition 1.1 in Chapter 5 of [4] show that we can use the divisors
$\Cal Z_b(f, \theta)$ to calculate the $\mu$-map. More precisely, since the
$f_i$
are disjoint, suppose that, for all choices of $b>0$
the intersection $$\Cal Z_{c-b}(f_{i_1}, \theta _{i_1}) \cap \dots \cap
\Cal Z_{c-b}(f_{i_{n-b}}, \theta _{i_{n-b}})= \emptyset,$$
that $\bigcap _{i=1}^n\Cal Z_c(f_i, \theta_i)=J$ is compact and
is contained in the Zariski open subset $M$ of $\frak M$, and that the $\Cal
Z_c(f_i, \theta _i)$ meet transversally along $J$. Here  by (2.6), $\frak
M_{c-b} = \emptyset$ if $b >(-p-3)/4$ and so it is enough to check the
above for
$0< b \leq (-p-3)/4$. In particular we must have $n>b$ for all $b\leq
(-p-3)/4$. Then we can define $\mu ([\Sigma])| J$ and arguments along the
lines of the proof of Theorem 1.12 in Chapter 5 of [4] show that  $$\gamma
_{w,p}(S)(f_1, \dots, f_n, [\Sigma], \dots, [\Sigma]) =  (\mu
([\Sigma])|J)^{d-n}.$$ Thus, Theorem 4.2 is a
consequence of the following:

\theorem{4.3} Suppose that $-p\geq 2(4p_g+2)$.  Let $f_1,
\dots, f_t$ denote distinct general fibers of $\pi$, and for each $i$ choose
$\theta _i$ a line bundle on $f_i$ with $\theta _i^2 = \scrO _{f_i}(\Delta)$.
Then: \roster
\item"{(i)}" For all $t\geq n$ and for all choices of $b>0$,  the
intersection
$$\Cal Z_{c-b}(f_{i_1}, \theta _{i_1}) \cap \dots \cap
\Cal Z_{c-b}(f_{i_{t-b}}, \theta _{i_{t-b}})= \emptyset.$$
\item"{(ii)}" If $t> n$, then $\bigcap _{i=1}^t\Cal Z_c(f_i, \theta_i)=
\emptyset$, and moreover $\bigcap _{i=1}^n\Cal Z_c(f_i, \theta_i)$ is compact
and is contained in the Zariski open subset $M$ of $\frak M$. The intersection
is
transverse and is a fiber $J$ of the map $M \to U\subseteq \Pee ^n$.
\item"{(iii)}" If $\Sigma$ is a smooth curve in $S$, then the restriction of
$\mu ([\Sigma])$ to each of the $m_2$ connected components  of $J$ is equal to
$(\Sigma \cdot f)\cdot [\Theta]$, where $\Theta$ is the theta divisor of the
component. \endroster
\endstatement
\medskip

We begin by determining the set-theoretic intersection of the $\Cal Z_c(f_i,
\theta _i)$. Recall that we have associated to $V$ the section $A$ and the
scheme $Z$ on the double cover $T$. The following is straightforward:

\lemma{4.4} Let $V\in \frak M$. Then $V\in \Cal Z_c(f, \theta)$ is and only if
either the section $A'$ of $\Bbb F$ corresponding to $V$ meets $B$
transversally at the point corresponding to $\theta$ or $\operatorname{Supp}Z
\cap \nu ^{-1}(f) \neq \emptyset$. \qed
\endstatement
\medskip

For the rest of the argument, we shall not worry about the case where the
sections pass through the exceptional point, since this can be handled by
symmetry. Arguing as in [4], Lemma 2.7 of Chapter 7, for a general choice of
fibers $f_1, \dots, f_t$ and line bundles $\theta _i$ on $f_i$, and for
all $s\leq r$, setting $H_i$ to be the hyperplane of sections in $|\sigma
+(N+s)f|$ passing through the point corresponding to $\theta _i$, then $H_1\cap
\dots \cap H_{N+2s +1} = \{A\}$ and the intersection of more than $N+2s +1$ of
the $H_i$ is empty. In particular, this means that the $f_i$ are chosen so that
the negative section $\sigma$ does not neet the branch divisor at any of the
$f_i$. Thus if $V$ is a bundle whose associated section  $A$ is
the negative section $\sigma$ of $\Bbb F_N$ and $V$ lies in the intersection
$\bigcap _{i=1}^t\Cal Z_c(f_i, \theta_i)$, then $\operatorname{Supp}Z$ meets
the
preimage in $T$ of $f_i$ for all $i$.

A counting argument as in the proof of  Proposition 2.9  of Chapter 7 in [4]
then establishes (i) and (ii) of Theorem 4.3, at least set-theoretically,
except
possibly for those $V$ such that the corresponding section of $\Bbb F_N$ is
the negative section. Likewise, arguments identical to the proof of (2.6) in
Chapter 7 of [4] show that the intersection of the divisor $\Cal Z_c(f_i,
\theta_i)$ with the Zariski open subset $M$ is reduced, so that the
intersection in (ii) of Theorem 4.3 is transverse. Thus the only new case to
consider is the possibility that the section associated to $V$ is the negative
section. We shall briefly outline the argument in this case; it is here that we
must assume that $-p$ is sufficiently large. In particular, recalling that $$n
=
\frac{d -p_g}2 = \frac{-p -4p_g-3}2,$$ and that we may assume that $b\leq
(-p-3)/4$, the condition $-p\geq 2(4p_g+2)$ insures that   $$\align n-b&\geq
n+\frac{p+3}4 = -\frac{p}4 -\frac{8p_g+3}4 \\ &\geq \frac{8p_g+4}2 -
\frac{8p_g+3}4 >0. \endalign$$ Thus the intersection in (i) is not over the
empty collection of divisors  $\Cal Z_{c-b}(f_i, \theta_i)$, i.e\. the
intersection is always contained in  a divisor  $\Cal Z_{c-b}(f_i, \theta_i)$.

Given $t$ general fibers $f_1, \dots, f_t$, we can assume that they do not
correspond to intersections of the negative section $\sigma$ with $B+\frak
d_1+\frak d_2$. Now let $V$ be a bundle whose associated section is $\sigma$.
Suppose further that $c_2(V) = c-b$ and that $V$ lies in the intersection
$\Cal Z_{c-b}(f_{i_1}, \theta _{i_1}) \cap \dots \cap
\Cal Z_{c-b}(f_{i_{t-b}}, \theta _{i_{t-b}})$ for $t\geq n$. We claim that, if
$0\leq b<(-p-3)/4$, then $-p \leq 4p_g +2N -2$. It clearly suffices to assume
that $t=n$. Let $\nu \:T \to S$ be the double cover corresponding to $\sigma$.
Then $\nu$ is not branched over the $f_i$. Writing $\nu ^*V$ as an extension
$$0
\to \scrO _T(D) \to \nu ^*V \to \scrO _T(\nu ^*\Delta - D)\otimes I_Z \to 0,$$
we
have by (2.5) that $-p =-p_1(\ad V) +4b = 4k-6N +1+ 2\ell (Z) + \delta+4b$.
Moreover $\operatorname{Supp}Z$ must meet each of the two components of $\nu
^{-1}(f_i)$ for $i=1, \dots, t$. Thus $\ell(Z) \geq 2n-2b$. But we also have
$$\align 2n &= d-p_g =-p - 4p_g -3\\ &= 4k-6N +1 + 2\ell (Z) + \delta -4(k-N-1)
-3+4b\\ &=-2N +2 + 2\ell (Z) +\delta+4b.
\endalign$$
Thus $2n \geq  -2N +2 + 4n-4b +\delta+4b$ and so $2n\leq 2N-2$. This says that
$$-p\leq 4p_g+2N-2.$$
On the other hand, $2(4p_g+2) =4p_g +4(k-N)\geq 4p_g +
4N>4p_g+2N-2$. Thus with our assumptions on $p$ there can be no bundle $V$ in
the intersection with an associated section equal to $\sigma$. Taking $b>0$
establishes (i) in the case where the section associated to $V$ is $\sigma$.
Taking $b=0$ shows that no bundle with associated section $\sigma$ lies in the
intersection $\bigcap _{i=1}^t\Cal Z_c(f_i, \theta _i)$ if $t\geq n$. This
establishes (i) and (ii).

Finally, to prove (iii) of Theorem 4.3, we must determine $\mu ([\Sigma])|J$,
where $J$ is a fiber of the map $M\to U$. The argument again closely parallels
that of [4] Chapter 7. A fiber of $M\to U$ determines and is determined by a
generic double cover $T\to S$. There is a divisor $D_0$ on $T$ such that
every $V$ in the fiber is of the form $\nu _*\scrO _T(D_0+F)\otimes \rho
^*\lambda$, for a line bundle $\lambda \in \Pic ^0C$. Fix a smooth holomorphic
multisection $\Sigma$ of $\pi$, transverse to the double cover $T\to S$. Let
$\Sigma'$ be the inverse image of $\Sigma$ under $\nu$. There is a commutative
diagram  $$\CD \Sigma ' @>{\nu}>> \Sigma\\
@V{\rho}VV   @VV{\pi}V\\
C @>{g}>> \Pee ^1.
\endCD$$
Let $\Cal P$ be the Poincar\'e bundle over $\operatorname{Pic}^0C\times C$.
Let $E$ be the  divisor  on $\Sigma '$ induced by $D_0 + F$. Then there is a
universal bundle over $\operatorname{Pic}^0C\times \Sigma$ of the form
$$(\Id \times \nu)_*\Bigl[\pi _2^*\scrO _{\Sigma '}(E)\otimes (\Id \times \rho
^*)\Cal P\Bigr].$$
Here $\pi _2\: \operatorname{Pic}^0C\times \Sigma' \to \Sigma '$ is the second
projection. The Chern classes of $\Cal V_E$ only depend on the numerical
equivalence class of $E$. Moreover, $p_1(\ad \Cal V_E) = p_1(\ad \Cal V_E
\otimes
q_2^*\lambda)$ for every line bundle $\lambda$ on $\Sigma$, where $q_2\:
\operatorname{Pic}^0C\times \Sigma \to \Sigma $ is the projection. This has the
effect of replacing $\scrO _{\Sigma '}(E)$ by $\scrO _{\Sigma '}(E)\otimes \nu
^*\lambda$. Replacing $\Sigma$ by $2\Sigma$, replaces $\deg E$ by $2\deg E$.
Thus we can assume that $\deg E$ is even, and then after twisting by an
appropriate $\lambda$ we can assume that $\deg E =0$. So as far as the Chern
classes are concerned we may as well assume that $E=0$, and we need to
calculate
$$p_1(\ad (\Id \times \nu)_*(\Id \times \rho)
^*\Cal P).$$

Now since cohomology commutes with flat base change, we have
$$(\Id \times \nu)_*(\Id \times \rho)
^*\Cal P= (\Id \times \pi)^*(\Id \times
g)_*\Cal P.$$
Thus we need to find
$$p_1(\ad (\Id \times \pi)^*(\Id \times
g)_*\Cal P) = (\Id \times \pi)^*p_1(\ad (\Id \times
g)_*\Cal P).$$
Let us first calculate $p_1(\ad (\Id \times g)_*\Cal P)$.
A straightforward calculation using e.g\. Lemma 2.14 of Chapter 7 of [4] shows
the following

\lemma{4.5} Let $f\:X \to Y$ be a double cover with $X$ and $Y$ smooth, and
let $\Upsilon$ be the line bundle on $Y$ defining the double cover, so that
$\Upsilon ^{\otimes 2} = \Cal {O}_Y(B)$, where $B$ is the branch locus of $f$.
If $D$ is a divisor on $X$, then
$$p_1(\ad f_*\scrO _X(D)) = c_1(\Upsilon )^2 - (f_*D)^2 + 2f_*(D^2).\qed$$
\endstatement
\medskip

Applying this in our situation, with $X = \Pic ^0C\times C$ and $f=\Id \times
g$, we see that $\Upsilon$ is the pullback of a divisor on $C$ and thus
$c_1(\Upsilon )^2 = 0$. So we are left with calculating (in the sense of
cycles)  $(\Id \times g)_*\Cal P$ and $(\Id \times g)_*[\Cal P]^2$. Also, as
in the proofs of (2.15) and (2.16) in Chapter 7 of [4],  $(\Id \times g)_*\Cal
P
=0$ and $[\Cal P ]^2 = -2r_1^*[\Theta] \cup r_2^*x$, where $x$ is the class of
a
point in $C$ and $\Theta$ is the theta divisor on $\Pic ^0C$, and $r_1, r_2$
are
the first and second projections on $\Pic ^0C\times C$. Thus
$$\align
2(\Id \times g)_*[\Cal P]^2 &=
2(\Id \times g)_*(-2r_1^*[\Theta] \cup r_2^*x)\\
&= -4p_1^*[\Theta] \cup p_2^*y,
\endalign$$
where $p_1$, $p_2$ are the first and second projections on $\Pic ^0C \times
\Pee^1$ and $y$ is the class of a point on $\Pee ^1$. It follows that
$-p_1(\ad \Cal V_E)/4 = q_1^*\Theta \cup q_2^*\pi ^*y$.

Hence the slant product of $-p_1(\ad \Cal V_E)/4$ with
$[\Sigma]$ is $(\deg \pi)\cdot[\Theta]$. Since $\deg \pi = (\Sigma \cdot f)$,
we have now established (iii) of Theorem 4.3. This concludes the proof of
Theorem
4.3 and thus of Theorem 4.1.  \qed


\Refs
\widestnumber\no{9}

\ref \no  1\by S. K. Donaldson and P. B. Kronheimer \book The Geometry of
Four-Manifolds \publ Clarendon \publaddr Oxford \yr 1990 \endref

\ref \no  2\by R. Friedman \paper Rank two vector bundles over regular
elliptic surfaces \jour Inventiones Math. \vol 96 \yr 1989 \pages 283--332
\endref

\ref \no 3\bysame \paper Vector bundles and $SO(3)$-invariants for
elliptic surfaces I \toappear \endref

\ref \no 4\by R. Friedman and J. W. Morgan \book Smooth 4-manifolds and complex
surfaces \toappear \endref


\endRefs
\enddocument